\documentclass[11pt]{article}

\usepackage{latexsym}
\usepackage{amssymb}

\newtheorem{theorem}{Theorem}

\newtheorem{example}{Example}

\begin{document}
{
\begin{center}
{\Large\bf
On some classical type Sobolev orthogonal polynomials.}
\end{center}
\begin{center}
{\bf S.M. Zagorodnyuk}
\end{center}
\section{Introduction.}

The theory of orthogonal polynomials on the real line has a lot of contributions and various applications, see~\cite{cit_50000_Gabor_Szego},
\cite{cit_98500_Freud_book},
\cite{cit_20000_Suetin}, \cite{cit_3000_Chihara}, \cite{cit_5200_Nevai_book},~\cite{cit_5000_Ismail}. 
One of its possible generalizations, the theory of Sobolev orthogonal polynomials, is nowadays of a considerable interest.
A brief account of the theory state, with references to several surveys on the theory, can be found in~\cite{cit_3500_CS_MB_JAT_2011},
see also~\cite{cit_5150_M_X}.
Consider the following second order differential equation:
\begin{equation}
\label{f1_10}
v (t) p_n''(t) + w (t) p_n'(t) + \lambda_n p_n(t) = 0,\qquad n\in\mathbb{Z}_+,
\end{equation}
where $v(t), w(t)$ are some real polynomials, $\lambda_n\in\mathbb{R}$. 
Orthogonal polynomials (OP) on the real line $p_n(t)$ which satisfy relation~(\ref{f1_10})
are called \textit{classical}. Namely, there are three such families: Jacobi OP, Hermite OP and Laguerre OP.
All of them are expressed in terms of the generalized hypergeometric function~\cite{cit_50_Bateman}.
The second-order operator in Equation~(\ref{f1_10}) can be replaced by a higher-order differential operator. The corresponding
polynomials were studied by Krall (Krall's polynomials) and many other mathematicians (see the book~\cite{cit_5100_Krall}, a
recent paper~\cite{cit_98150_DI} and references therein).

There already appeared Sobolev orthogonal polynomials satisfying differential equations, see~\cite{cit_300_AAMP} for the case of Jacobi-Sobolev OP.
However, the theory of classical Sobolev OP seems to be at its beginning, and it is not complete as for Equation~(\ref{f1_10}).
Let $\rho$ and $\xi$ be some fixed non-negative integers. Denote
\begin{equation}
\label{f1_20}
L_\rho y(t) = \sum_{j=0}^\rho l_j(t) y^{(j)}(t),\quad D_\xi y(t) = \sum_{k=0}^\xi d_k(t) y^{(k)}(t),
\end{equation}
where $l_j,d_k$ and $y$ are real polynomials of $t$: $l_\rho\not= 0$, $d_\xi\not= 0$.
Thus, $L_\rho$ and $D_\xi$ are linear differential operators with real polynomial coefficients having orders $\rho$ and $\xi$,
respectively.
Consider the following differential equation:
\begin{equation}
\label{f1_30}
L_\rho y(t) + \lambda D_\xi y(t) = 0,
\end{equation}
where $\lambda\in\mathbb{R}$. 
Polynomial solutions of equation~(\ref{f1_30}), with $\rho=\xi=2$, were studied by 
Sawyer and Chaundy in~\cite{cit_15000_Sawyer},\cite{cit_1000_Chaundy}.
Chaundy also trusted that the method and arguments in~\cite{cit_1000_Chaundy} could be enlarged to deal with equations of
higher order.
In~\cite{cit_97100_Z} there appeared orthogonal polynomials which are solutions of a differential equation of the
form~(\ref{f1_30}) with $\rho= 4$, $\xi=2$, and with $\rho= 3$, $\xi=1$.

In this paper we present a way to construct new families of Sobolev orthogonal polynomials using 
known families of orthogonal polynomials on the real line.
In particular, we consider two families of hypergeometric polynomials:
\begin{equation}
\label{f1_40}
\mathcal{L}_n(x) = \mathcal{L}_n(x;q,r) = {}_2 F_2(-n,1;q,r;x),\qquad q,r>0,\ n\in\mathbb{Z}_+;
\end{equation}
and
\begin{equation}
\label{f1_50}
\mathcal{P}_n(x) = \mathcal{P}_n(x;a,b,c) = {}_3 F_2(-n,n-1+a+b,1;a,c;x),\quad a,b,c>0,\ n\in\mathbb{Z}_+.
\end{equation}
Polynomials $\mathcal{L}_n, \mathcal{P}_n$
generalize Laguerre and Jacobi polynomials, respectively.
These polynomials satisfy differential
equations of the form~(\ref{f1_30}) with $\rho=3,\xi=1$. 
For positive integer values of the parameters $r,c$, polynomials $\mathcal{L}_n, \mathcal{P}_n$
are Sobolev orthogonal polynomials
with some explicitly given matrix weights. We study some basic properties of these polynomials, which follow from
their hypergeometric nature. In particular, some linear
recurrence relations are obtained by Fasenmeier's method.
A perfect exposition of the latter method can be found in Rainville's book~\cite{cit_5150_Rainville}.
The situation in our case is complicated by the fact that
polynomials depend on several positive parameters. Therefore coefficients of the corresponding linear algebraic systems
(for unknown coefficients of a recurrence relation) depend on these parameters. Thus, the leading coefficients in the
Gauss elimination method can become unavailable for analysis.
Finally, we discuss some possible generalizations of our constructions.

\noindent
{\bf Notations. }
As usual, we denote by $\mathbb{R}, \mathbb{C}, \mathbb{N}, \mathbb{Z}, \mathbb{Z}_+$,
the sets of real numbers, complex numbers, positive integers, integers and non-negative integers,
respectively. 
By $\mathbb{Z}_{k,l}$ we mean all integers $j$ satisfying the following inequality:
$k\leq j\leq l$; ($k,l\in\mathbb{Z}$).
By $\mathbb{P}$ we denote the set of all polynomials with complex coefficients.
For a complex number $c$ we denote
$(c)_0 = 1$, $(c)_1=c$, $(c)_k = c(c+1)...(c+k-1)$, $k\in\mathbb{N}$ (\textit{the shifted factorial}).
The generalized hypergeometric function is denoted by
$$ {}_m F_n(a_1,...,a_m; b_1,...,b_n;x) = \sum_{k=0}^\infty \frac{(a_1)_k ... (a_m)_k}{(b_1)_k ... (b_n)_k} \frac{x^k}{k!}, $$
where $m,n\in\mathbb{N}$, $a_j,b_l\in\mathbb{C}$.
By $\Gamma(z)$ and $\mathrm{B}(z)$ we denote the gamma function and the beta function, respectively.

\section{Hypergeometric Sobolev orthogonal polynomials.}

Let $\{ p_n(x) \}_{n=0}^\infty$ ($p_n$ has degree $n$ and real coefficients) be orthogonal polynomials
on $[a,b]\subseteq \mathbb{R}$ with respect to a weight function $w(x) (\geq 0)$:
\begin{equation}
\label{f2_10}
\int_a^b p_n(x) p_m(x) w(x) dx = A_n \delta_{n,m},\qquad A_n>0,\quad n,m\in\mathbb{Z}_+.
\end{equation}
We do not need for $p_n$ to have positive leading coefficients. The weight $w$ is assumed to be continuous on $(a,b)$.

Fix a positive integer $\xi$. Consider the following differential equation:
\begin{equation}
\label{f2_25}
D_\xi y(t) = p_n(t),
\end{equation}
where $D_\xi$ is defined as in Equation~(\ref{f1_20}), and $n\in\mathbb{Z}_+$. The following assumption will play a key role in
what follows.

\noindent
\textbf{Condition 1.}
\textit{Suppose that for each $n\in\mathbb{Z}_+$, the differential equation~(\ref{f2_25}) has a real $n$-th degree polynomial
solution $y(t) = y_n(t)$.}

If Condition~1 is satisfied, by relations~(\ref{f2_10}),(\ref{f2_25}) we immediately obtain that 
$\{ y_n(x) \}_{n=0}^\infty$ are Sobolev orthogonal polynomials:

$$ \int_a^b (y_n(x), y_n'(x),\ldots, y_n^{(\xi)}(x)) M(x) 
\left(
\begin{array}{cccc}
y_m(x) \\
y_m'(x) \\
\vdots \\
y_m^{(\xi)}(x)
\end{array}
\right)
w(x) dx = A_n \delta_{n,m}, $$
\begin{equation}
\label{f2_30}
n,m\in\mathbb{Z}_+,
\end{equation}
where
\begin{equation}
\label{f2_35}
M(x) := 
\left(
\begin{array}{cccc}
d_0(x) \\
d_1(x) \\
\vdots \\
d_{\xi}(x)
\end{array}
\right)
(d_0(x), d_1(x),\ldots, d_\xi(x)),\qquad x\in (a,b).
\end{equation}
Fix an arbitrary $\rho\in\mathbb{N}$, and consider the differential operator $L_\rho$ as in~(\ref{f1_20}).
If $\{ p_n(x) \}_{n=0}^\infty$ satisfy the differential equation:
\begin{equation}
\label{f2_45}
L_\rho p_n(t) = \lambda_n p_n(t),\qquad \lambda_n>0,\quad n\in\mathbb{Z_+},
\end{equation}
then $y_n(t)$ satisfy the following differential equation:
\begin{equation}
\label{f2_50}
L_\rho D_\xi y_n(t) = \lambda_n D_\xi y_n(t),\qquad  n\in\mathbb{Z_+}.
\end{equation}
In this case $\{ y_n(x) \}_{n=0}^\infty$ are classical type Sobolev orthogonal polynomials.

In order to construct the above classical type polynomials, it is necessary to find differential equations~(\ref{f2_25})
which satisfy Condition~1. The following two theorems ensure that such equations indeed exist.
Let
$$ P_n^{(\alpha,\beta)} (t) = 
\left(
\begin{array}{cc} n+\alpha \\
n 
\end{array}
\right)
{}_2 F_1\left(-n,n+\alpha+\beta + 1; \alpha+1; \frac{1-t}{2} \right), $$
$$ \alpha,\beta > -1,\ n\in\mathbb{Z}_+, $$
denote the Jacobi polynomials, and
$$ L_n^\alpha (x) = 
\left(
\begin{array}{cc} n+\alpha \\
n 
\end{array}
\right)
{}_1 F_1(-n;\alpha + 1;x),\qquad \alpha>-1,\ n\in\mathbb{Z}_+, $$
denote the (generalized) Laguerre polynomials (\cite{cit_50_Bateman}).

\begin{theorem}
\label{2_1}
Let $r$ be an arbitrary positive integer greater than $1$.
Polynomials $y_n(x) = \mathcal{L}_n(x;q,r)$ ($q>0$) satisfy orthogonality relations~(\ref{f2_30}),(\ref{f2_35})
with $\xi = r-1$; $a=0$; $b=\infty$; $w(x) = x^{q-1} e^{-x}$,
\begin{equation}
\label{f2_55}
A_n = 
\frac{ \left( (r-1)! \right)^2 \Gamma(q+n) }
{
\left\{
\left(
\begin{array}{cc} n+q-1 \\
n 
\end{array}
\right)
\right\}^2
n!
},\qquad n\in\mathbb{Z}_+;
\end{equation}
\begin{equation}
\label{f2_60}
d_k(x) = 
\left(
\frac{ (r-1)! }{ k! }
\right)^2
\frac{1}{(r-1-k)!}
x^k,\qquad k\in\mathbb{Z}_{0,r-1}.
\end{equation}
Moreover, polynomials $y_n$ satisfy differential equation~(\ref{f2_50}), where
$L_\rho$ and $D_\xi$ are defined as in Relation~(\ref{f1_20}),
with $\rho = 2$; $l_2(t)=t$; $l_1(t)=q-t$; $l_0(t)=0$; $\lambda_n = -n$. 
\end{theorem}
\textbf{Proof.}
Observe that
$$ \left( x^{r-1} \mathcal{L}_n(x;q,r) \right)^{(r-1)} = 
\left( 
\sum_{k=0}^\infty \frac{ (-n)_k }{ (q)_k (r)_k } x^{k+r-1}
\right)^{(r-1)} = 
$$
$$ =
(r-1)! \sum_{k=0}^\infty \frac{ (-n)_k }{ (q)_k } \frac{ x^k }{ k! } 
=
\frac{ (r-1)! }{ \left(
\begin{array}{cc} n+q-1 \\
n 
\end{array}
\right)
} 
L_n^{q-1}(x).
$$
Thus, in our case we can choose
$$ p_n(x) = \frac{ (r-1)! }{ \left(
\begin{array}{cc} n+q-1 \\
n 
\end{array}
\right)
} 
L_n^{q-1}(x). $$
Relation~(\ref{f2_25}) holds with
\begin{equation}
\label{f2_62}
D_\xi y(t) = \left( x^{r-1} y(t) \right)^{(r-1)} = \sum_{j=0}^{r-1} \left(
\frac{ (r-1)! }{ j! }
\right)^2
\frac{1}{ (r-1-j)! } x^j y^{(j)}(t). 
\end{equation}
Therefore Condition~1 is satisfied.
Moreover, since $L_n^{q-1}(x)$ satisfy the differential equation
$$ xy'' + (q-x) y' = -ny, $$
then relation~(\ref{f2_45}) holds with
$\rho = 2$, $l_2(x)=x$, $l_1(x)=q-x$, $l_0(x)=0$, and $\lambda_n = -n$.
Using the orthogonality conditions for the Laguerre polynomials~(\cite{cit_50_Bateman}) and 
relations~(\ref{f2_30}),(\ref{f2_50}) we complete the proof.
$\Box$

%In a similar manner we obtain the following theorem.

\begin{theorem}
\label{2_2}
Let $c$ be an arbitrary positive integer greater than $1$.
Polynomials $y_n(x) = \mathcal{P}_n(x;a,b,c)$ ($a,b>0$) satisfy orthogonality relations~(\ref{f2_30}),(\ref{f2_35})
with $\xi = c-1$, $a=0,b=1$, $w(x) =  (2x)^{a-1} (2-2x)^{b-1}$,
$$ A_n = 
\frac{ \left( (c-1)! \right)^2 2^{a+b-2} \Gamma(a+n) \Gamma(b+n) }
{
\left\{
\left(
\begin{array}{cc} n+a-1 \\
n 
\end{array}
\right)
\right\}^2
(2n+a+b-1) n! \Gamma(n+a+b-1)
},\ n\in\mathbb{N}, $$
$$ 
A_0 = \left( (c-1)! \right)^2 2^{a+b-2} \mathrm{B}(a,b); $$
\begin{equation}
\label{f2_70}
d_k(x) = 
\left(
\frac{ (c-1)! }{ k! }
\right)^2
\frac{1}{(c-1-k)!}
x^k,\qquad k\in\mathbb{Z}_{0,c-1}.
\end{equation}
Polynomials $y_n$ satisfy differential equation~(\ref{f2_50}), where
$L_\rho$ and $D_\xi$ are defined as in Relation~(\ref{f1_20}),
with $\rho = 2$, $l_2(t)=t(1-t)$, $l_1(t)= a-t(a+b)$, $l_0(t)=0$, $\lambda_n = -n(n+a+b-1)$. 
\end{theorem}
\textbf{Proof.}
The operator $D_\xi$ will be as in Relation~(\ref{f2_62}) but with $c$ instead of $r$; $\xi=c-1$. Observe that
$$ \left( x^{c-1} \mathcal{P}_n(x;a,b,c) \right)^{(c-1)}  
= (c-1)! \sum_{k=0}^\infty \frac{ (-n)_k (n-1+a+b)_k }{ (a)_k } \frac{ x^k }{ k! } = $$
$$ = (c-1)! {}_2 F_1 (-n,n-1+a+b;a;x) = 
\frac{ (c-1)! }{ \left(
\begin{array}{cc} n+a-1 \\
n 
\end{array}
\right)
} 
P_n^{(a-1,b-1)}(1-2x).
$$
Thus, we can choose
$$ p_n(x)  
=
\frac{ (c-1)! }{ \left(
\begin{array}{cc} n+a-1 \\
n 
\end{array}
\right)
} 
P_n^{(a-1,b-1)}(1-2x).
$$
Notice that polynomials
$$ \widehat{P}_0^{(\alpha,\beta)} (x) := \frac{1}{ \sqrt{ 2^{\alpha+\beta+1} \mathrm{B}(\alpha+1,\beta+1) } }, $$
$$ \widehat{P}_n^{(\alpha,\beta)} (x) := \sqrt{ \frac{ (2n+\alpha+\beta+1) n! \Gamma(n+\alpha+\beta+1) }{ 2^{\alpha+\beta+1}
\Gamma(n+\alpha+1) \Gamma(n+\beta+1) } } P_n^{(\alpha,\beta)} (x),\quad n\in\mathbb{N}, $$
are orthonormal on $[-1,1]$ (having positive leading coefficients) with respect to the weight
$\mathbf{w}(x)=(1-x)^\alpha (1+x)^\beta$, $\alpha,\beta>-1$.
Therefore
$$ \mathbf{P}_0^{(\alpha,\beta)} (t) := (-1)^n \sqrt{2} \frac{1}{ \sqrt{ 2^{\alpha+\beta+1} B(\alpha+1,\beta+1) } }, $$
$$ \mathbf{P}_n^{(\alpha,\beta)} (t) := (-1)^n \sqrt{2}
\sqrt{ \frac{ (2n+\alpha+\beta+1) n! \Gamma(n+\alpha+\beta+1) }{ 2^{\alpha+\beta+1}
\Gamma(n+\alpha+1) \Gamma(n+\beta+1) } } P_n^{(\alpha,\beta)} (1-2t), $$
$$ n\in\mathbb{N}, $$
are orthonormal polynomials on $[0,1]$ with respect to $\mathbf{w}(-2t+1) = (2t)^\alpha (2-2t)^\beta$,
see~\cite[p. 29]{cit_50000_Gabor_Szego}. 
Thus, 
$p_n(x)$ are orthogonal polynomials on $[0,1]$ with the weight $w(x)$ and constants $A_n$, as in the statement of the theorem.
Consequently, Condition~1 is satisfied.
Since polynomials $P_n^{(\alpha,\beta)}(x)$ satisfy the differential equation:
$$ (1-x^2) y'' + [\beta-\alpha-(\alpha+\beta+2)x] y' + n(n+\alpha+\beta+1) y = 0, $$
then
$p_n(t)$ satisfy the following differential equation:
$$ t(1-t) w'' + [a - (a+b)t] w' = - n(n+a+b-1) w. $$
Therefore differential equation~(\ref{f2_50}) holds.
$\Box$

In the preceding theorems we have seen that the order of the corresponding differential equation~(\ref{f2_50}) for
$\mathcal{L}_n(x;q,r)$ and $\mathcal{P}_n(x;a,b,c)$ depends on their parameters.
It turns out that polynomials $\mathcal{L}_n(x;q,r)$ and $\mathcal{P}_n(x;a,b,c)$ satisfy differential equations
of the form~(\ref{f1_30}) with $\rho=3$, $\xi=1$.
Moreover the latter equations exist for arbitrary parameters of these polynomials.

\begin{theorem}
\label{2_3}
Polynomials $\mathcal{L}_n(x;q,r)$ ($q,r>0$) satisfy the following differential equation:
\begin{equation}
\label{f2_75}
x^2 y''' + (q+r+1-x) xy'' + (qr-2x) y' + n \left\{
xy' + y
\right\} = 0.
\end{equation}
Polynomials $\mathcal{P}_n(x;a,b,c)$ ($a,b,c>0$) satisfy the following differential equation:
$$ (1-x)x^2 y''' + (a+c+1 - (a+b+3)x) xy'' + (ac - 2(a+b)x) y' + $$
\begin{equation}
\label{f2_80}
+ n(n+a+b-1) \left\{
xy' + y
\right\} = 0.
\end{equation}
\end{theorem}
\textbf{Proof.}
Since the generalized hypergeometric function satisfies a linear differential equation of a special form~\cite{cit_50_Bateman},
we obtain that 
polynomials $\mathcal{L}_n(x;q,r)$ ($q,r>0$) satisfy the following differential equation:
\begin{equation}
\label{f2_85}
\left\{
\delta (\delta + q -1) (\delta + r -1) - x(\delta - n)(\delta + 1)
\right\} 
y(x) = 0,
\end{equation}
while polynomials $\mathcal{P}_n(x;a,b,c)$ ($a,b,c>0$) satisfy the differential equation
\begin{equation}
\label{f2_90}
\left\{
\delta (\delta + a -1) (\delta + c -1) - x(\delta - n)(\delta+n-1+a+b)(\delta + 1)
\right\}
y(x) = 0,
\end{equation}
where $\delta = \frac{d}{dx}$.
By gathering the corresponding derivatives and simplifications we get the desired result.
$\Box$

The hypergeometric nature of polynomials $\mathcal{L}_n,\mathcal{P}_n$ allows to give an integral representation for them.
\begin{theorem}
\label{2_4}
Polynomials $\mathcal{L}_n(z;q,r)$, with $q>0$ and $r>1$, admit the following integral representation:
\begin{equation}
\label{f2_93}
\mathcal{L}_n(z;q,r) = 
\frac{ (r-1) }{ 
\left(
\begin{array}{cc} n+q-1\\
n\end{array}
\right)
}
\int_0^1
(1-t)^{r-2}
L_n^{q-1} (zt) dt,\qquad \forall z\in\mathbb{C}. 
\end{equation}
Polynomials $\mathcal{P}_n(x;a,b,c)$, with $a,b>0$ and $c>1$, have the following representation:
$$ \mathcal{P}_n(z;a,b,c) = 
\frac{ (c-1) }{ 
\left(
\begin{array}{cc} n+a-1\\
n\end{array}
\right)
}
\int_0^1
(1-t)^{c-2}
P_n^{(a-1,b-1)} (1-2zt) dt, $$
\begin{equation}
\label{f2_95}
\forall z\in\mathbb{C}:\ |z|<1. 
\end{equation}
\end{theorem}
\textbf{Proof.}
Let $q>0,r>1$, and $a,b>0$, $c>1$.
By Theorem~28 in~\cite[p. 93]{cit_5150_Rainville} we may write:
$$ {}_2 F_2(1,-n;r,q;z) = (r-1)\int_0^1 (1-t)^{r-2} {}_1 F_1(-n;q;zt) dt,\quad z\in\mathbb{C}; $$
$$ {}_3 F_2(1,-n,n-1+a+b;c,a;z) = $$
$$ = (c-1)\int_0^1 (1-t)^{c-2} {}_2 F_1(-n,n-1+a+b;a;zt) dt,\quad |z|<1. $$
It remains to use the hypergeometric representations for the corresponding polynomials to get
representations~(\ref{f2_93}),(\ref{f2_95}).
$\Box$

The following example shows that polynomials $\mathcal{L}_n,\mathcal{P}_n$ need not to be orthogonal on the real line.

\begin{example}
\label{e2_1}
Observe that
$$ \mathcal{L}_2(z;q,r) = 1 - \frac{2}{qr} z + \frac{2}{q(q+1) r(r+1)} z^2. $$
The discriminant of the following quadratic equation
$$ 2 z^2 - 2(q+1)(r+1) z + q(q+1)r(r+1) = 0, $$
is equal to
\begin{equation}
\label{f2_97}
D_L(q,r) := 4(q+1)(r+1) \left\{
-qr + q + r + 1
\right\}.
\end{equation}
We have $D_L(r,r) = 4(r+1)^2 (-r^2 + 2r + 1)\rightarrow -\infty$, as $r\rightarrow\infty$.
In particular, $D_L(3,3) = -2 < 0$, $D_L(1+\sqrt{2},1+\sqrt{2})=0$. Hence,
$\mathcal{L}_2(z;3,3)$ has complex roots, while $\mathcal{L}_2(z;1+\sqrt{2},1+\sqrt{2})$ has a multiple root.

Notice that
$$ \mathcal{P}_2(z;a,b,c) = 1 - \frac{ 2(a+b+1) }{ ac } z + \frac{ 2(a+b+1)(a+b+2) }{a(a+1) c(c+1)} z^2. $$
The discriminant of the following quadratic equation
$$ 2(a+b+1)(a+b+2) z^2 - 2 (a+b+1) (a+1)(c+1) z + a(a+1)c(c+1) = 0, $$
is equal to
\begin{equation}
\label{f2_102}
D_P(a,b,c) := 4(a+b+1)(a+1)(c+1) (
-a^2 c -2ac + 2a - abc + a^2 + ab + bc + b + c + 1 )
\end{equation}
We have 
$$ D_P(a,a,a) = 4(2a+1)(a+1)^2 (-2a^3 + a^2 + 4a +1)\rightarrow -\infty,\qquad \mbox{as } a\rightarrow\infty. $$
In particular, $D_P(3,3,3) < 0$.
\end{example}

We shall now obtain recurrence relations for polynomials $\mathcal{P}_n(z;a,b,c)$, $\mathcal{L}_n(x;q,r)$.
Let $\delta$ be a positive integer. The above polynomials admit the following generalizations:
$$ \mathbf{L}_n(x) = \mathbf{L}_n(x;q,r_1,...,r_\delta) = {}_{\delta+1} F_{\delta+1} (-n,1,...,1;q,r_1,...,r_\delta;x),  $$ 
\begin{equation}
\label{f2_103}
q,r_1,...,r_\delta>0,\ n\in\mathbb{Z}_+;
\end{equation}
and
$$ \mathbf{P}_n(x) = \mathbf{P}_n(x;a,b,c_1,...,c_\delta) =  $$
\begin{equation}
\label{f2_105}
{}_{\delta+2} F_{\delta+1}(-n,n-1+a+b,1,...,1;a,c_1,...,c_\delta;x),\qquad a,b,c_1,...,c_\delta >0,\ n\in\mathbb{Z}_+.
\end{equation}

Fix an arbitrary integer $n\geq 2$.
Observe that
\begin{equation}
\label{f2_107}
\mathbf{P}_{n+1}(x) = \sum_{k=0}^\infty (-n-1)_k (n+a+b)_k \frac{ (k!)^{\delta -1} }{ (a)_k (c_1)_k ... (c_\delta)_k } x^k = 
\sum_{k=0}^\infty \varepsilon_k(n), 
\end{equation}
where 
$$ \varepsilon_k(n) := (-n-1)_k (n+a+b)_k \frac{ (k!)^{\delta -1} }{ (a)_k (c_1)_k ... (c_\delta)_k } x^k. $$
According to Fasenmeier's method (\cite{cit_5150_Rainville}) we need to express
$\mathbf{P}_{n}(x)$, $\mathbf{P}_{n-1}(x)$, $\mathbf{P}_{n-2}(x)$, $x \mathbf{P}_{n}(x)$ and $x \mathbf{P}_{n-1}(x)$
in terms of $\varepsilon_k(n)$.
Since
$$ (-n)_k = (-n-1)_k \frac{(n+1-k)}{(n+1)}, $$
$$ (n+a+b-1)_k = (n+a+b)_k \frac{(n+a+b-1)}{(n+a+b+k-1)},\qquad k\in\mathbb{Z}_+, $$
we obtain that
\begin{equation}
\label{f2_109}
\mathbf{P}_{n}(x) = \sum_{k=0}^\infty \varepsilon_k(n) \frac{(n+1-k)}{(n+1)} \frac{(n+a+b-1)}{(n+a+b-1+k)}. 
\end{equation}
In a similar way we get
$$ \mathbf{P}_{n-1}(x) = $$
\begin{equation}
\label{f2_112}
= \sum_{k=0}^\infty \varepsilon_k(n) \frac{(n+1-k)(n-k)}{(n+1)n} \frac{(n+a+b-2)(n+a+b-1)}{(n+a+b-2+k)(n+a+b-1+k)}, 
\end{equation}
$$ \mathbf{P}_{n-2}(x) = \sum_{k=0}^\infty \varepsilon_k(n) \frac{(n+1-k)(n-k)(n-1-k)}{(n+1)n(n-1)} * $$
\begin{equation}
\label{f2_115}
* \frac{(n+a+b-3)(n+a+b-2)(n+a+b-1)}
{(n+a+b-3+k)(n+a+b-2+k)(n+a+b-1+k)}. 
\end{equation}
Notice that
\begin{equation}
\label{f2_117}
(r)_{k-1} = \frac{1}{(r+k-1)} (r)_k,\qquad k\in\mathbb{N}.
\end{equation}
Then
$$ x \mathbf{P}_{n}(x) = \sum_{k=0}^\infty (-n)_k (n+a+b-1)_k \frac{ (k!)^{\delta -1} }{ (a)_k (c_1)_k ... (c_\delta)_k } x^{k+1} = $$
$$ = \sum_{j=1}^\infty (-n)_{j-1} (n+a+b-1)_{j-1} \frac{ ((j-1)!)^{\delta -1} }{ (a)_{j-1} (c_1)_{j-1} ... (c_\delta)_{j-1} } x^j = $$
$$ = \sum_{k=1}^\infty \varepsilon_k(n) (-1) k \frac{ (n+a+b-1) }{ (n+1)(n+a+b-2+k)(n+a+b-1+k) } * $$
\begin{equation}
\label{f2_119}
* \frac{ (a+k-1)(c_1+k-1)...(c_\delta+k-1) }{ k^\delta }. 
\end{equation}
By similar arguments we obtain that
$$ x \mathbf{P}_{n-1}(x) = \sum_{k=1}^\infty \varepsilon_k(n) (-1) k * $$
$$ * \frac{ (n+1-k)(n+a+b-2)(n+a+b-1) }{ (n+1)n (n+a+b-3+k)(n+a+b-2+k)(n+a+b-1+k) } * $$
\begin{equation}
\label{f2_122}
* \frac{ (a+k-1)(c_1+k-1)...(c_\delta+k-1) }{ k^\delta }. 
\end{equation}

Consider the following expression:
$$ R_1(x) := 
\varphi_1(n) \mathbf{P}_{n-2}(x) + \varphi_2(n) \mathbf{P}_{n-1}(x) + \varphi_3(n) \mathbf{P}_{n}(x) + \varphi_4(n) \mathbf{P}_{n+1}(x) + $$
$$ + \varphi_5(n) x \mathbf{P}_{n-1}(x) + \varphi_6(n) x \mathbf{P}_{n}(x), $$
where $\varphi_j = \varphi_j(n)$ are free real parameters.
We shall try to select these parameters in such a way that ensures $R_1(x)=0$ for all $x$.
By relations~(\ref{f2_107}),(\ref{f2_109}),(\ref{f2_112}),(\ref{f2_115}) and (\ref{f2_119}),(\ref{f2_122}) the expression
$R_1(x)$ takes the following form:
$$ R_1(x) = 
\varepsilon_0(n) \left(
\varphi_1(n) + \varphi_2(n) + \varphi_3(n) + \varphi_4(n)
\right) + $$
$$ + \sum_{k=1}^\infty
\varepsilon_k(n)
\left\{
\varphi_1(n) \frac{ (n+1-k)(n-k)(n-1-k) }{ (n+1)n(n-1) } * \right. $$
$$ * \frac{ (n+a+b-3)(n+a+b-2)(n+a+b-1) }{ (n+a+b-3+k)(n+a+b-2+k)(n+a+b-1+k) }
+ $$
$$ + \varphi_2(n) \frac{ (n+1-k)(n-k) (n+a+b-2)(n+a+b-1) }{ (n+1)n (n+a+b-2+k)(n+a+b-1+k) } + $$
$$ + \varphi_3(n) \frac{ (n+1-k) (n+a+b-1) }{ (n+1) (n+a+b-1+k) } + \varphi_4(n)
+
$$ 
$$ + \left(
\varphi_5(n) \frac{ (n+1-k)(n+a+b-2) }{ n(n+a+b-3+k) } + \varphi_6(n)
\right) * $$
$$ * (-1) \frac{ (n+a+b-1) }{ (n+1)(n+a+b-2+k)(n+a+b-1+k) } * $$
\begin{equation}
\label{f2_125}
\left. * \frac{ (a+k-1)(c_1+k-1)...(c_\delta+k-1) }{ k^{\delta-1} }\right\}.
\end{equation}
Denote the expression in the last brackets $\{ ... \}$ by $E_k(n)$.
Observe that $\varepsilon_k(n) = 0$, for all $k\geq n+2$.
If
\begin{equation}
\label{f2_127}
\varphi_1(n) + \varphi_2(n) + \varphi_3(n) + \varphi_4(n) = 0,
\end{equation}
and
\begin{equation}
\label{f2_129}
E_k(n) = 0,\qquad k\in\mathbb{Z}_{1,n+1},
\end{equation}
then $R_1(x)\equiv 0$.
As it is common when applying Fasenmeier's method, $E_k(n)$ is a rational function of $k$:
$E_k(n) = \frac{P_k(n)}{Q_k(n)}$, where $P,Q$ are polynomials of $k$.
Thus, instead of~(\ref{f2_129}) it is enough to check
the equality $P_k(n)=0$ for a fixed number $N$ ($N>\deg P_k(n)$) of distinct points $k$.

The problem in our case is that $\deg P_k(n)$ depends on $\delta$. Therefore, in general we shall obtain a large
system of linear algebraic equations with coefficients depending on parameters $a,b,c_1,...,c_\delta$.
Besides huge expressions, it is not clear how to guarantee that the leading coefficients during the Gauss elimination
will be nonzero.
Moreover we only have $6$ unknowns. 
We can conjecture that a recurrence relation for $\mathbf{P}_n$ should include $3+\delta$ subsequent polynomials.

In what follows we assume that $\delta = 1$, $c_1 = c\in\mathbb{R}$.
Multiply~(\ref{f2_129}) by
$(n+a+b-3+k)(n+a+b-2+k)(n+a+b-1+k)$, and denote
$$ \psi_1(n) = \frac{ (n+a+b-3)(n+a+b-2)(n+a+b-1) }{ (n+1)n(n-1) } \varphi_1(n), $$
$$ \psi_2(n) = \frac{ (n+a+b-2)(n+a+b-1) }{ (n+1)n } \varphi_2(n), $$
$$ \psi_3(n) = \frac{ (n+a+b-1) }{ (n+1) } \varphi_2(n),\quad \psi_4(n) = \varphi_4(n), $$ 
$$ \psi_5(n) = \frac{ (-1) (n+a+b-2)(n+a+b-1) }{ (n+1)n } \varphi_5(n), $$
$$ \psi_6(n) = \frac{ (-1) (n+a+b-1) }{ (n+1) } \varphi_6(n), $$
to get
$$ \psi_1(n) (n+1-k)(n-k)(n-1-k) + \psi_2(n) (n+1-k)(n-k)(n+a+b-3+k) + $$
$$ + \psi_3(n) (n+1-k)(n+a+b-3+k)(n+a+b-2+k) + $$
$$ + \psi_4(n) (n+a+b-3+k)(n+a+b-2+k)(n+a+b-1+k) + $$
$$ + \psi_5(n) (n+1-k)(a+k-1)(c+k-1) + $$
\begin{equation}
\label{f2_131}
+ \psi_6(n) (n+a+b-3+k)(a+k-1)(c+k-1) = 0,\qquad k\in\mathbb{Z}_{1,n+1}. 
\end{equation}
The left side of relation~(\ref{f2_131}) is a polynomial of degree $\leq 3$ of $k$. Thus, if~(\ref{f2_131}) 
holds for $4$ distinct values of $k$, then it holds for all complex $k$.
These values should be selected carefully to obtain as simple equations as possible. 
We choose $k$ to be $n+1,n,n-1$ and $k=0$. The value $k=n-2$ would lead to a more complicated equation.
After the substitution of these values and some simplifications, involving the application of~(\ref{f2_127}), we obtain:
\begin{equation}
\label{f2_132}
\psi_4(n) (2n+a+b-1)(2n+a+b) + \psi_6(n) (a+n)(c+n) = 0,
\end{equation}
$$ \psi_3(n) (2n+a+b-3)(2n+a+b-2) + $$
$$ + \psi_4(n) (2n+a+b-3)(2n+a+b-2)(2n+a+b-1) + $$
\begin{equation}
\label{f2_134}
+ \psi_5(n) (a+n-1)(c+n-1) + \psi_6(n) (2n+a+b-3) (a+n-1)(c+n-1) = 0,
\end{equation}
$$ \psi_2(n) 2 (2n+a+b-4) + \psi_3(n) 2 (2n+a+b-4)(2n+a+b-3) + $$
$$ + \psi_4(n) (2n+a+b-4)(2n+a+b-3)(2n+a+b-2) + $$
\begin{equation}
\label{f2_136}
+ \psi_5(n) 2 (a+n-2)(c+n-2) + \psi_6(n) (2n+a+b-4) (a+n-2)(c+n-2) = 0,
\end{equation}
\begin{equation}
\label{f2_138}
\psi_5(n) (n+1) (a-1)(c-1) + \psi_6(n) (n+a+b-3) (a-1)(c-1) = 0.
\end{equation}
We can express $\psi_4$ and $\psi_5$ in terms of $\psi_6$, by equations~(\ref{f2_132}) and~(\ref{f2_138}).
Then $\psi_3$ is expressed in terms of $\psi_6$ by equation~(\ref{f2_134}),
and $\psi_2$ is expressed in terms of $\psi_6$ by equation~(\ref{f2_136}).
Thus, we can calculate $\varphi_2,\varphi_3,\varphi_4,\varphi_5$ by an arbitrary $\varphi_6$.
Finally, $\varphi_1$ is determined by relation~(\ref{f2_127}).
We choose
\begin{equation}
\label{f2_145}
\varphi_6(n) = (2n+a+b-1)(2n+a+b)(n+1)(n+a+b-2). 
\end{equation}
Then
\begin{equation}
\label{f2_147}
\varphi_4(n) = (a+n)(c+n) (n+a+b-1)(n+a+b-2); 
\end{equation}
\begin{equation}
\label{f2_149}
\varphi_5(n) = (2n+a+b-1)(2n+a+b)(n+1)(n+a+b-2); 
\end{equation}
$$ \varphi_3(n) = (-1) (2n+a+b-1) (n+a+b-2) 
\left\{
(n+1)(a+n)(c+n) + 
\right.
$$
$$ + (n+a+b-3) (2n+a+b) \frac{ (a+n-1)(c+n-1) }{ (2n+a+b-3)(2n+a+b-2) } - $$
\begin{equation}
\label{f2_151}
\left. 
- (2n+a+b) (n+1) \frac{ (a+n-1)(c+n-1) }{ (2n+a+b-2) }
\right\}.
\end{equation}
Denote by $I_1$ the expression standing in brackets $\{ ... \}$ in~(\ref{f2_151}).
By writing $(2n+a+b) = (2n+a+b-3) + 3$, and $(2n+a+b) = (2n+a+b-2) + 2$, in the last two summands of $I_1$,
we can simplify the fractions to get
$$ I_1 = (n+1)(2n+a+c-1) + (-n+a+b-2) \frac{ (a+n-1)(c+n-1) }{ (2n+a+b-2) } - $$
\begin{equation}
\label{f2_153}
- 3n \frac{ (a+n-1)(c+n-1) }{ (2n+a+b-3)(2n+a+b-2) }.
\end{equation}
Writing $(-n+a+b-2) = (2n+a+b-2) - 3n$, in the second summand
we simplify the fraction to get
\begin{equation}
\label{f2_155}
I_1 = n (2n+a+c-1) + (a+n)(c+n) - 3n \frac{ (a+n-1)(c+n-1) }{ (2n+a+b-3) }.
\end{equation}
Therefore
$$ \varphi_3(n) = (-1) (2n+a+b-1) (n+a+b-2) * $$
\begin{equation}
\label{f2_157}
* \left\{ n (2n+a+c-1) + (a+n)(c+n) - 3n \frac{ (a+n-1)(c+n-1) }{ (2n+a+b-3) } \right\}.
\end{equation}
For $\varphi_2$ we obtain the following expression:
$$ \varphi_2(n) = n [ (2n+a+b-3)(2n+a+b-1) ( n(2n+a+c-1) + (a+n)(c+n) ) - $$
$$ - 3n (2n+a+b-1) (a+n-1)(c+n-1) - $$
$$ - \frac{1}{2}
(a+n)(c+n)(2n+a+b-3)(2n+a+b-2)(n+1) - $$
$$ - (n+a+b-3)(2n+a+b-1)(2n+a+b) \frac{ (a+n-2)(c+n-2) }{ (2n+a+b-4) } + $$
\begin{equation}
\label{f2_159}
+ \frac{1}{2}
(2n+a+b-1)(2n+a+b)(n+1)(a+n-2)(c+n-2) ].
\end{equation}
Finally, we get
\begin{equation}
\label{f2_161}
\varphi_1(n) = - \varphi_2(n) - \varphi_3(n) - \varphi_4(n).
\end{equation}
We shall not substitute for $\varphi_1,\varphi_2,\varphi_3$ in the last expression, since the resulting expression is very large.

\begin{theorem}
\label{2_5}
Polynomials $\mathcal{P}_n(x) = \mathcal{P}_n(x;a,b,c)$ ($a,b,c>0$) satisfy the following recurrence relation:
$$ \varphi_1(n) \mathcal{P}_{n-2}(x) + \varphi_2(n) \mathcal{P}_{n-1}(x) + \varphi_3(n) \mathcal{P}_{n}(x) + \varphi_4(n) \mathcal{P}_{n+1}(x) + $$
\begin{equation}
\label{f2_163}
+ \varphi_5(n) x \mathcal{P}_{n-1}(x) + \varphi_6(n) x \mathcal{P}_{n}(x),\qquad n\in\mathbb{Z}_+, 
\end{equation}
where $\varphi_j(n)$ ($j\in\mathbb{Z}_{1,6}$) are defined by 
relations~(\ref{f2_145}),(\ref{f2_149}),(\ref{f2_147}),(\ref{f2_157}),(\ref{f2_159}),(\ref{f2_161}) for
$n\geq 2$. For $n=0,1$ the coefficients $\varphi_1, \varphi_4,\varphi_5,\varphi_6$ are defined by the same formulas.
Moreover, $\varphi_2(0) = 0$, $\varphi_3(0) = - ac (a+b-1)(a+b-2)$;
$$ \varphi_2(1) = (a+b-1)(a+b+1) (a+c+1 +(a+1)(c+1)) - $$ 
$$ -3ac(a+b+1) - (a+1)(c+1)(a+b-1)(a+b), $$
$$ \varphi_3(1) = - (a+b+1)(a+b-1) (a+c+1 + (a+1)(c+1)) + 3ac (a+b+1), $$
and $\mathcal{P}_{-2}(x) = \mathcal{P}_{-1}(x) = 0$.
\end{theorem}
\textbf{Proof.} For $n\geq 2$ formula~(\ref{f2_163}) follows from the preceding considerations. For $n=0,1$
it can be verified directly. $\Box$

Fortunately, in order to obtain a recurrence relation for $\mathcal{L}_n(x;q,r)$ ($q,r>0$) we need not proceed in the same way.
We shall make use of the following property:
\begin{equation}
\label{f2_167}
\mathcal{L}_n(x;q,r) = \lim\limits_{b\to +\infty} \mathcal{P}_n \left(
\frac{x}{b};q,b,r
\right),\qquad n\in\mathbb{Z}_+,\ x\in\mathbb{R};\ q,r>0.
\end{equation}
Relation~(\ref{f2_167}) readily follows from the hypergeometric representations of the corresponding polynomials.

\begin{theorem}
\label{2_6}
Polynomials $\mathcal{L}_n(x) = \mathcal{L}_n(x;q,r)$ ($q,r>0$) satisfy the following recurrence relation:
$$ (n-1)n \mathcal{L}_{n-2}(x) - n(3n+q+r-2) \mathcal{L}_{n-1}(x) + $$
$$ + (n(2n+q+r-1) + (n+q)(n+r)) \mathcal{L}_n(x) - (n+q)(n+r) \mathcal{L}_{n+1}(x) + $$
\begin{equation}
\label{f2_305}
+ n x \mathcal{L}_{n-1}(x) - (n+1) x \mathcal{L}_n(x) = 0,\qquad n\in\mathbb{Z}_+,
\end{equation}
where $\mathcal{L}_{-2}(x) = \mathcal{L}_{-1}(x)=0$.
\end{theorem}
\textbf{Proof.}
Write relation~(\ref{f2_163}) for
$\mathcal{P}_n \left(
\frac{x}{b};q,b,r
\right)$, then
divide it by $b^3$ and pass to the limit as $b\rightarrow +\infty$.
$\Box$

As we have seen, polynomials $\mathcal{P}_n(x;a,b,c)$ and $\mathcal{L}_n(x;q,r)$ are (generalized) eigenvectors of pencils of differential operators,
as well as of pencils of difference operators, for all possible parameters. However, the orthogonality relations
were proved for positive integer values of parameters only. Thus, there appears the following problem.

\noindent
\textbf{Open problem~1}. Do there exist Sobolev type orthogonality relations for $\mathcal{P}_n(x;a,b,c)$ and $\mathcal{L}_n(x;q,r)$
for all positive values of $a,b,c$ and $q,r$?

\textbf{Generalizations.} Let us return to polynomials 
$\mathbf{L}_n(x) = \mathbf{L}_n(x;q,r_1,...,r_\delta)$, and $\mathbf{P}_n(x) = \mathbf{P}_n(x;a,b,c_1,...,c_\delta)$,
($q,r_1,...,r_\delta; a,b,c_1,...,c_\delta >0$, $\delta\in\mathbb{N}$) defined in~(\ref{f2_103}),(\ref{f2_105}).
If parameters $r_1,...,r_\delta$ and $c_1,...,c_\delta$ are positive integers, then Condition~1 is satisfied for these polynomials. 
In fact, the following expression:
$$ \left( x^{r_\delta-1} \mathbf{L}_n(x;q,r_1,...,r_\delta) \right)^{(r_\delta-1)}, $$
is, up to a constant factor, the polynomial $\mathbf{L}_n(x;q,r_1,...,r_{\delta-1})$.
Thus, after applying this operation $\delta$ times we come to the Laguerre polynomials.
The case of $\mathbf{P}_n(x)$ is similar. 
Differential equations can be also written for these polynomials, for arbitrary positive parameters. 
As for recurrence relations, we have seen the difficulties which
arise here.

In the general case, one can write equation~(\ref{f2_25}) with an unknown polynomial $y(t)$ and equate the coefficients
by the same powers of $x$ (the same simple but powerful idea was used in~\cite{cit_5_Azad}).
The problem here is to get a compact representation of the corresponding solutions.
This will be studied elsewhere.

In Equation~(\ref{f2_25}) instead of orthogonal polynomials on the real line one can consider
(bi)orthogonal rational functions $p_n(t)$ (\cite{cit_98000_Zhedanov_JAT}),
orthogonal polynomials on the unit circle, etc.

It is also of interest to study the zero distribution of polynomials $\mathcal{L}_n, \mathcal{P}_n$.
Zeros of hypergeometric polynomials are intensively studied nowadays, see, e.g.~\cite{cit_98500_ZSW_Proc_AMS} \cite{cit_500_BM}.

\begin{center}
{\large\bf 
On some classical type Sobolev orthogonal polynomials.}
\end{center}
\begin{center}
{\bf S.M. Zagorodnyuk}
\end{center}

In this paper we propose a way to construct classical type Sobolev orthogonal polynomials.
We consider two families of hypergeometric polynomials:
${}_2 F_2(-n,1;q,r;x)$ and
${}_3 F_2(-n,n-1+a+b,1;a,c;x)$
($a,b,c,q,r>0$, $n=0,1,...$), which generalize Laguerre and Jacobi polynomials, respectively.
These polynomials satisfy higher-order differential
equations of the following form: $L y + \lambda_n D y = 0$, where $L,D$ are linear
differential operators with polynomial coefficients not depending on $n$.
For positive integer values of the parameters $r,c$ these polynomials are Sobolev orthogonal polynomials
with some explicitly given measures. Some basic properties of these polynomials, including recurrence relations,
are obtained.

}


\begin{thebibliography}{99}

\bibitem{cit_300_AAMP}
Arves\'u, J.; \'Alvarez-Nodarse, R.; Marcell\'an, F.; Pan, K. 
Jacobi-Sobolev-type orthogonal polynomials: second-order differential equation and zeros. 
J. Comput. Appl. Math. 90 (1998), no. 2, 135--156. 

\bibitem{cit_5_Azad}
Azad H., Laradji A., Mustafa M. T. Polynomial solutions of differential equations. Adv. Difference Equ. 2011:58 (2011), 12 pp.

\bibitem{cit_50_Bateman}
Erd\'elyi, Arthur; Magnus, Wilhelm; Oberhettinger, Fritz; Tricomi, Francesco G. 
Higher transcendental functions. Vols. I, II. Based, in part, on notes left by Harry Bateman. 
McGraw-Hill Book Company, Inc., New York-Toronto-London, 1953. {\rm xxvi}+302, {\rm xvii}+396 pp.

\bibitem{cit_500_BM}
Bracciali, Cleonice F.; Moreno-Balc\'azar, Juan Jos\'e. 
On the zeros of a class of generalized hypergeometric polynomials. Appl. Math. Comput. 253 (2015), 151--158.

\bibitem{cit_1000_Chaundy}
Chaundy, T. W. Second-order linear differential equations with polynomial solutions. Quart. J. Math., Oxford Ser. (2) 4, (1953). 81--95.


\bibitem{cit_3000_Chihara}
Chihara, T. S. An introduction to orthogonal polynomials. Mathematics and its Applications, Vol. 13. 
Gordon and Breach Science Publishers, New York-London-Paris, 1978. xii+249 pp.

\bibitem{cit_3500_CS_MB_JAT_2011}
Costas-Santos, R. S.; Moreno-Balc\'azar, J. J. The semiclassical Sobolev orthogonal polynomials: a general approach. 
J. Approx. Theory 163 (2011), no. 1, 65--83.


\bibitem{cit_98150_DI}
Dur\'an, Antonio J.; de la Iglesia, Manuel D. 
Differential equations for discrete Jacobi-Sobolev orthogonal polynomials. J. Spectr. Theory 8 (2018), no. 1, 191--234


\bibitem{cit_98500_Freud_book}
Freud, G\'eza. Orthogonal polynomials.  
Pergamon Press, Oxford New York Toronto Sydney Braunschweig, 1971.


\bibitem{cit_5000_Ismail}
Ismail, Mourad E. H. Classical and quantum orthogonal polynomials in one variable. With two chapters by Walter Van Assche. 
With a foreword by Richard A. Askey. Encyclopedia of Mathematics and its Applications, 98. Cambridge University Press, Cambridge, 2005. xviii+706 pp.

\bibitem{cit_5100_Krall}
Krall, Allan M. Hilbert space, boundary value problems and orthogonal polynomials. 
Operator Theory: Advances and Applications, 133. Birkhäuser Verlag, Basel, 2002. xiv+352 pp.

\bibitem{cit_5150_M_X}
Marcell\'an, Francisco; Xu, Yuan. On Sobolev orthogonal polynomials. Expo. Math. 33 (2015), no. 3, 308--352.

\bibitem{cit_5200_Nevai_book}
Nevai, Paul G. Orthogonal polynomials.  
Memoirs of the AMS, Number~213, 1979, 185~pp.


\bibitem{cit_5150_Rainville}
Rainville, Earl D. Special functions. Reprint of 1960 first edition. Chelsea Publishing Co., Bronx, N.Y., 1971. {\rm xii}+365 pp.


\bibitem{cit_15000_Sawyer}
Sawyer, W. W. Differential equations with polynomial solutions. Quart. J. Math., Oxford Ser. 20, (1949). 22--30.


\bibitem{cit_20000_Suetin}
Suetin, P. K. Classical orthogonal polynomials.
Third edition. Fizmatlit, Moscow, 2005. 480 pp. (Russian)


\bibitem{cit_50000_Gabor_Szego}
Szeg\"o, G\'abor. Orthogonal polynomials. Fourth edition. 
American Mathematical Society, Colloquium Publications, Vol. XXIII. American Mathematical Society, Providence, R.I., 1975. xiii+432 pp.

\bibitem{cit_97100_Z}
Zagorodnyuk S. M., Difference equations related to Jacobi-type pencils.---
{\it J. Difference Equ. Appl.}, {\bf 24}, no. 10 (2018), 1664--1684.

\bibitem{cit_98000_Zhedanov_JAT}
Zhedanov A., Biorthogonal rational functions and the generalized eigenvalue problem.
J. Approx. Theory, 101 (1999), no. 2, pp. 303--329.


\bibitem{cit_98500_ZSW_Proc_AMS}
Zhou, Jian-Rong; Srivastava, H. M.; Wang, Zhi-Gang. 
Asymptotic distributions of the zeros of a family of hypergeometric polynomials. Proc. Amer. Math. Soc. 140 (2012), no. 7, 2333--2346.

\end{thebibliography}
\end{document}